%%%%%%%%%%%%%%%%%%%%%%%%%%%%%%%%%%%%%%%%%%%%%%%%%%%%%%%%%%%%%%%%%%%%%%%%%%%%% 
%%%%%%%%%%%%%%%%%%%%%%%%%%VOREINSTELLUNGEN%%%%%%%%%%%%%%%%%%%%%%%%%%%%%%%%%%% 
%%%%%%%%%%%%%%%%%%%%%%%%%%%%%%%%%%%%%%%%%%%%%%%%%%%%%%%%%%%%%%%%%%%%%%%%%%%%% 
\input amstex 
\documentstyle{amsppt} 
\hoffset=.2true in 
\hsize= 6 true in  
\vsize= 9 true in 
\vcorrection{0.2in} 
\NoRunningHeads 
%%%%%%%%%%%%%%%%%%%%%%%%%%%%%%%%%%%%%%%%%%%%%%%%%%%%%%%%%%%%%%%%%%%%%%%%%%%%%% 

\def\ni{\noindent}% 
\def\ts{\thinspace}% 
\def\CC{{\Bbb C}}% 
\def\PP{{\Bbb P}}% 
\def\NN{{\Bbb N}}% 
\def\QQ{{\Bbb Q}}%
\def\ZZ{{\Bbb Z}}% 
\def\"#1{{\accent"7F #1\penalty10000\hskip 0pt plus 0pt}} % Umlaute! %
 
\scrollmode 
%%%%%%%%%%%%%%%%%%%%%%%%%%%%%%%%%%%%%%%%%%%%%%%%%%%%%%%%%%%%%%%%%%%%%%%%%%%%%% 
\topmatter  
\title Geometric Diffeomorphism Finiteness in Low Dimensions 
and Homotopy Group Finiteness
%\vskip2mm
\endtitle 

\author Wilderich Tuschmann%%\footnote"$^\bigstar$"
%{With thanks
%to the Max-Planck Institute Leipzig
%or its support and hospitality. \hfill{$\,$}}
%\vskip2mm
\endauthor
 
\address Max-Planck-Institute for Mathematics in the Sciences, 
Inselstrasse, D-04103 Leipzig, Germany 
\endaddress 
\email tusch\@mis.mpg.de 
\endemail 

%\vskip3mm

\date August 1999
\enddate

\abstract  
Our main result asserts that 
for any given numbers $C$ and $D$ 
the class of simply connected closed smooth manifolds of dimension 
$m<7$
which admit a Riemannian metric 
with sectional curvature bounded in absolute value by $\vert K \vert\le C$ 
and diameter uniformly bounded from above by $D$ contains 
only  finitely many diffeomorphism types.
Thus in these dimensions the lower positive bound on volume
in Cheeger's Finiteness Theorem
can be replaced by a purely topological condition, 
simply-connectedness.
In dimension $4$ instead of simply-connectedness
here only non-vanishing of the Euler characteristic has to be required.

As a topological corollary we obtain that for $k+l<7$
there are
over a given smooth closed $l$-manifold 
only finitely many principal $T^k$ bundles
with simply connected and non-diffeomorphic total spaces.

Furthermore, for any given numbers $C$ 
and $D$ and any dimension $m$ it is shown
that for each $i\in\NN$ there are up to isomorphism
always only finitely many possibilities for the $i$th homotopy group
of a simply connected closed $m$-manifold which admits a metric
with curvature $\vert K \vert\le C$ and diameter $\le D$.
\endabstract       

\endtopmatter         
\document 
 
\baselineskip=.16in

\head  Introduction 
\endhead 
 
\vskip3mm

\ni 
The main result of this note asserts that in dimensions less than $7$
the presence of a metric
with given curvature and diameter bounds
suffices to restrict the diffeomorphism type of a 
simply connected closed manifold
to finitely many possibilities:

\vskip4mm 
 
\ni 
{\bf Theorem 1.\ts}  
{\sl 
For any given numbers $C$ and $D$
there is only a finite number of diffeomorphism types of
simply connected closed smooth
$m$-manifolds, $m< 7$,
which admit Riemannian 
metrics with sectional curvature  
$|K|\le C$ and diameter $\le D$.
}

\vskip4mm

Note that Theorem $1$ does not require a lower positive bound
on volume or injectivity radius. 

\vskip3mm

The above result implies by Myer's theorem that
given $m<7$ and any $\delta>0$,
there is only a finite number of diffeomorphism types of  
simply connected closed smooth
$m$-dimen\-sio\-nal manifolds $M$
which admit Riemannian metrics with Ricci curvature $Ricc\ge\delta>0$
and sectional curvature $K\le 1$.

This explains in particular
why $7$ is the first dimension where infinite
sequences of closed simply connected manifolds 
of mutually distinct diffeomorphism type and
uniformly positively pinched sectional curvature 
(see [AW], [E]) can appear.

%\vskip4mm
 
In view of Theorem $1$ it is also interesting to note 
that there are simply connected four-manifolds
of fixed homeomorphism type 
which admit infinitely many distinct smooth structures 
(for example $\CC\PP^2 \# \ts 9 \overline{\CC\PP^2}$, see [FM]).
(In dimension $4$ Theorem $1$ actually holds
for closed smooth manifolds with non-zero Euler characteristic, see below.)

\vskip3mm

In the relevant dimensions
Theorem $1$ complements and improves other geometric finiteness results.
It shows that for dimensions $m<7$
in Cheeger's Finiteness Theorem the assumption of a lower positive
bound on volume can be replaced by a topological condition,
namely, simply-connectedness.
Theorem $1$ also shows that for these dimensions in the
$\pi_2$-Finiteness Theorem from [PT] 
the requirement that the second homotopy group be finite is actually
not needed.

%\vskip4mm

\vfill\eject

The finiteness theorems in Riemannian geometry
which require at most bounds on volume, curvature, and diameter
can be stated as follows: 
For manifolds $M$ of a given dimension $m$, the conditions

\vskip2mm

\item{$\bullet$} $vol(M)\ge v>0$, $|K(M)|\le C$ and diam$(M)\le D$
imply finiteness of diffeomorphism types (Cheeger ([C])
and Peters ([Pt]));
this conclusion moreover continues to hold under the conditions
$vol(M)\ge v>0$, $\int_M\vert R\vert \sp {m/2}\leq C$, $|Ricc_M|\le C'$,
diam$(M)\le D$ (Anderson and Cheeger ([AC]));

\item{$\bullet$} $vol(M)\ge v>0$, $K(M)\ge C$ diam$(M)\le D$ 
imply finiteness of homotopy types (Grove and Petersen ([GP])),
homeomorphism types (Perelman ([Pr]) and 
Lipschitz homeomorphism types (Perelman, unpublished);
if in addition $m\ge 4$, these conditions imply finiteness of homeomorphism,
and, if $m>4$, finiteness of diffeomorphism types
(Grove, Petersen and Wu ([GPW]));

\item{$\bullet$} $\pi_2(M)$ finite, $\pi_1(M)=0$, 
$|K(M)|\le C$ and diam$(M)\le D$ 
imply finiteness of diffeomorphism types ([PT]);

\item{$\bullet$} $K(M)\ge C$ and diam$(M)\le D$ 
imply a uniform bound for the total Betti number (Gromov ([G]).

\vskip4mm

Note that in every dimension $m\ge 7$ there are
counterexamples to Theorem $1$:

%\vskip2mm 

Uniformly positively pinched sequences of mutually non-diffeomorphic
Aloff-Wallach
or Eschenburg spaces (see [AW], [E])
and their products with spheres of appropiate dimension show that
the theorem does not hold in dimension $7$, nor in any dimension $m\ge 9$.

To obtain counterexamples to the validity of Theorem $1$ in dimension $8$
one can proceed as follows: 

Starting with the connected sum
$M:=S^3 \times \CC\PP^2 \# \ts S^3 \times \CC\PP^2$, using the Gysin sequence
one sees that for any pair of relatively prime integers $p,q\in\ZZ$
there is a circle bundle $S^1\to E_{p,q} \to M$ over $M$ whose
total space is simply connected and whose fourth cohomology is
isomorphic to $\ZZ_p \times  \ZZ_q$. In particular infinitely
many non-diffeomorphic total spaces arise. 
It only remains to observe
that given any Riemannian metric on $M$,
one can easily construct on each total space $E_{p,q}$ a metric with similar
curvature and diameter bounds as the metric on $M$.

\vskip3mm 

\ni
{\bf Question.\ts} When can the upper curvature bound in Theorem $1$
be discarded ?

\vskip2mm

In dimension $3$ the answer might well depend on the
still unresolved Poincar\' e Conjecture.

Combining Gromov's Betti number theorem and Freedman's classification
of simply connected topological four-manifolds (see [FQ]) one can see
that Theorem $1$ holds in dimension $4$ under the conditions 
$K\ge C$ and diam $\le D$, provided that
one replaces diffeomorphism by homeomorphism types.
(The reader might wish to compare this with the finiteness theorem of 
Grove-Petersen-Wu).

In dimension $5$ the answer to the above question seems to be positive.
Using
Barden's diffeomorphism classification of simply connected five-manifolds
([B]) and again the Betti number theorem,
it can at least be shown that
in dimension $5$ the conclusion of Theorem $1$ holds
under the conditions
$K(M)\ge C$, diam$(M)\le D$, and $\vert$Tor $H_2(M;\ZZ)\vert\le C'$.

In the six-dimensional case the upper curvature bound
is absolutely necessary for Theorem $1$ to hold.
This follows from the fact that
Karsten Grove and Wolfgang Ziller constructed nonnegatively curved metrics
on an infinite sequence of 
non-diffeomorphic $S^2$ bundles over $S^4$ (see [GZ]).

\vskip3mm

The proof of Theorem $1$ uses in dimensions $5$ and $6$
the classification results of Barden ([B]) and \v Zubr ([Z]).
Its geometric ingredients consists in collapsing arguments and the
following homotopy group finiteness theorem
which, in contrast to Theorem $1$, is valid in any dimension:

\vskip4mm

\ni 
{\bf Theorem 2.\ts} 
{\sl Given $m\in\NN$, $C$ and $D$, for each natural number $i\ge 2$
there exists a finite set $\Pi_i=\Pi_i(m,C,D)$ of
isomorphism classes of
finitely generated Abelian groups that satisfy:

If $M$ is a simply connected closed smooth
$m$-dimen\-sio\-nal manifold which admits a Riemannian 
metric with sectional curvature  
$|K(M)|\le C$ and diameter diam$(M)\le D$, then $\pi_i(M)\in\Pi_i$.
}

\vskip2mm 

Theorem $2$ follows from
a geometric classification theorem
for simply connected closed manifolds of [PT] (see below)
and improves
Theorem $0.3$ in [R]
from finiteness of possibilities
for the rational homotopy groups $\pi_i(M)\otimes\QQ$ 
to finiteness of possibilities for the homotopy groups $\pi_i(M)$ themselves.

\vskip4mm

I would like to mention that Fang and Rong announced
an independent proof of Theorem $1$,
and it is my pleasure to thank Anton Petrunin and Matthias Kreck,
without whom I also would never have found the reference [H],
for help with topology.

\vfill\eject 
 
\head  Proofs
\endhead 
 
\vskip4mm
\ni 
{\bf Proof of Theorem $1$.\ts}
In dimension one and two the theorem is true by trivial reasons. 
Since a simply connected closed three-manifold is a homotopy three-sphere,
in dimension three the $\pi_2$-Finiteness Theorem of [PT] applies. 
Thus we are left
with the cases where the dimension $m$ is equal to $4$, $5$, or $6$. 

\vskip2mm

First note that in a given dimension $m$ Theorem $1$ is true
for all closed smooth $m$-manifolds with non-zero Euler characteristic.
This can be seen as follows:

Suppose that
there exists a sequence $(M_n)_{n\in\NN}$ of
pairwise non-diffeomorphic closed smooth Riemannian $m$-manifolds
with uniformly bounded curvatures and diameters
and non-zero Euler characteristic.
Cheeger's Finiteness Theorem
implies that this sequence must collapse, i.e., it must hold that
$vol(M_n)\to 0$ as $n\to\infty$.
In particular, by [CG] for $n$ sufficiently large all $M_n$
admit a pure $F$-structure of positive rank and thus have vanishing
Euler characteristic, which is a contradiction.

\vskip2mm

Since the Euler characteristic of a 
closed $4$-manifold with finite fundamental group is at least $2$,
Theorem $1$ therefore holds in particular
for closed smooth $4$-manifolds with arbitrarily large,
but finite order of the fundamental group.

\vskip2mm

In the five-dimensional case Barden's classification
([B]) says that
simply connected $5$-manifolds $M$ are classified up to diffeomorphism
by the second homology group $H_2(M;\ZZ)$ and an invariant $i(M)$ which is 
obtained as follows: Regarding the second Stiefel-Whitney class of $M$
as a homomorphism $w:H_2(M)\to\ZZ_2$, one may arrange $w$ to be non-zero
on at most one element of a certain "basis" of $H_2$. This element has order
$2^i$ for some $i$, and this $i$ is the invariant $i(M)$.

Since there are for a given finitely generated 
%% Abelian %%
group $H$ always only
finitely many homomorphisms to 
%% $\ZZ_2$ %%
a fixed finite group, 
it follows in particular
that there are always only finitely many distinct diffeomorphism
types of closed smooth simply connected $5$-manifolds 
with a given second homology group.
This observation, combined with Theorem $2$ 
(and the Hurewicz theorem) proves Theorem $1$ for $m=5$. 

\vskip2mm

To prove Theorem $1$ in dimension six,
let us first note that
\v Zubr's diffeomorphism classification 
of closed oriented simply connected six-manifolds ([Z]) implies that
a given class of closed smooth simply connected manifolds
of dimension $6$ will contain at most finitely many diffeomorphism types
if it satisfies the following conditions:

For the members $M$ of this class there are, up to isomophism,
only finitely many possibilities for the second homology group
$H_2(M;\ZZ)$, 
the third Betti number $b_3(M;\ZZ)$,
the cup form $\mu_M$ (a symmetric trilinear form
$H^2(M;\ZZ)\otimes H^2(M;\ZZ)\otimes H^2(M;\ZZ)\to H^6(M;\ZZ)\cong \ZZ$ given
by the cup product evaluated on the orientation class,
which determines the multiplicative structure of the cohomology ring),
the first Pontryagin class $p_1(M)$ (which here is integral,
i.e., $p_1(M)\in H^4(M;\ZZ)$).

Now suppose that for some $C$ and $D$ there
exists an infinite sequence of pairwise non-diffeomorphic closed smooth
simply connected Riemannian $6$-manifolds $M_n$
with curvature $\vert K \vert \le C$
and diameter $\le D$.
By what has been said above, 
this sequence must collapse, and 
we may suppose that all manifolds
$M_n$ carry a pure $F$-structure of positive rank and
have vanishing Euler characteristic.

Since $M_n$ is simply connected,
by [CG] (compare [PT], [R])
the $F$-structure on
each $M_n$ is given by an effective smooth torus action without fixed points.
Since all orbits of this action have positive dimension,
it is easy to see that each $M_n$ thus also admits a
fixed-point free {\it circle} action.

Now ([H], Lemma 3.2) implies that each $M_n$ has vanishing
trilinear cup form $\mu$ and that the first Pontryagin class of $M_n$
is torsion, i.e., $p_1(M_n)\in$ Tor $ H^4(M_n;\ZZ)$.
By Poincar\'e duality, Tor $ H^4(M_n;\ZZ)$ $\cong$ Tor $H_2(M_n;\ZZ)$.
Also note that $\chi(M_n)=0$ implies that $b_3(M_n)=2\ts b_2(M_n)+2$.

But by Theorem $2$ there are only finitely many possibilities for the second
homology group of all $M_n$.
Combining this with the above facts it follows
that our collapsing sequence contains at most finitely non-diffeomorphic
manifolds, which yields the desired contradiction.
$\qed$

\vskip4mm

Since the total space of a principal
torus bundle over a Riemannian manifold
always carries a metric with similar curvature and diameter bounds as the base,
Theorem $1$ has the following
\vskip3mm
\noindent
{\bf Topological Corollary.\ts}
For $k+l<7$
there are
over a given smooth closed $l$-manifold 
only finitely many principal $T^k$ bundles
with simply connected and non-diffeomorphic total spaces.

\vfill\eject

Theorem $2$ is a direct corollary of the
following classification theorem
for simply connected closed manifolds from [PT]:

\vskip4mm

\ni 
{\bf Theorem ([PT]).\ts}
{\sl For given $m\in\NN$, $C$ and $D$, 
there exists a finite number of closed smooth manifolds $E_l$
such that any simply connected closed  
$m$-dimen\-sio\-nal manifold $M$ admitting a Riemannian 
metric with sectional curvature  
$|K|\le C$ and diameter $\le D$ is diffeomorphic to a  
factor space $M=E_l/T^{k_l}$, where $0\le k_l=$ dim $E_l-m$  
and $T^{k_l}$ acts freely on $E_l$. 
}

\vskip4mm

\ni 
{\bf Proof of Theorem 2.\ts}
Fix numbers $m$, $C$, and $D$, 
and let $M$ be a closed smooth simply connected of dimension $m$ which
admits a Riemannian metric with sectional curvature  
$|K|\le C$ and diameter $\le D$.
By the above theorem
there exists a closed smooth simply connected manifold $E$ 
such that $E$ is diffeomorphic to
the total space of a principal $T^k$ bundle over $M$.
(In the terminology of [PT],
the manifold $E$ is the {\it universal torus bundle}
of $M$, where $\pi_2(E)$ is finite and $0\le k=b_2(M)$.

Since tori are aspherical, the homotopy exact sequence 
$$\cdots \to \pi_i(T^k)\to\pi_i(E)\to\pi_i(M)\to\cdots
\to 0=\pi_2(T^k)\to\pi_2(E)\to\pi_2(M)\to
\pi_1(T^k)\cong\ZZ^k\to 0$$
of the principal bundle $T^k\to E\to M$ then shows
that for $3\le i\in \NN$ the homotopy group $\pi_i(M)$ is isomorphic
to $\pi_i(E)$, and that $\pi_2(M)$ is isomorphic to $\pi_2(E)\oplus\ZZ^k$.

Since by the above theorem for given numbers $m$, $C$, and $D$
there are only finitely many non-diffeomorphic manifolds $E$,
in noting that the homotopy groups of a closed simply connected manifold
are finitely generated one sees that the proof of
Theorem $2$ is complete.
$\qed$

\enddocument